# NANYANG RESEARCH PROGRAMME

# Unbalancing Lights or Gale-Berlekamp Switch Game

Le Viet Hung

Raffles Institution

Xu Yu

Victoria Junior College

Assistant Professor Kiah Han Mao

Nanyang Technological University

### *Abstract*

The Gale-Berlekamp Light Switching Game is played on a square board of lights. Each light has two states, either on or off. There is a switch to every row and column. Turning this switch would change the state of all the lights on that row or column. The aim of the game is to minimise the imbalance in the board, defined to be the absolute difference between the number of lights switched on and that of lights switched off. We investigate variants of the game for an *m x n* matrix with *n* even and *m ≤ n*. We provide a constructive proof that for any *m x n* rectangle matrix *A*, there exists $x \in (\pm 1)^n$ and $y \in (\pm 1)^m$ such that $|yAx| \le 2$. i.e. column and row switches to reduce the imbalance to at most *2*, construct a complete Python routine to find these switches, and test run the algorithm against randomly generated initial board configurations. We then expand the game to a three-dimensional *n x n x n* cube, with corresponding row, column, and layer switches. We define a minimum threshold $P_n$, such that the imbalance can always be reduced to at most $P_n$, for all initial states of the cube. We then provide an existential proof that $P_2 = 2$ and $P_4 = 4$.

## 1 INTRODUCTION

The Gale-Berlekamp Light Switching game is a well-known mathematical game proposed by American mathematician Elwyn Berlekamp. There are two main directions for the study: either maximising or minimising the imbalance in the board given an initial configuration of the board. For decades, mathematicians have studied the game in the context of square boards, but few have looked at the cases for rectangle and cubes. Inspired by the proof of J. Beck and J. Spencer on minimising unbalanced lights in the game [1], we aim to examine the variants of Gale-Berlekamp Light Switching Game, including both rectangle boards and simple cube configurations, and attempt to minimise the imbalance in each case.

The Gale-Berlekamp Light Switching game is played on a square board of lights. Each light has two states, either on or off. There is a switch to every row and column. Turning this switch would change the state of all the lights on that row or column. The aim of the game is to minimise the imbalance in the board. The **imbalance** is defined as the absolute difference between the number of lights switched on and that of lights switched off.

For an *n x n* square board, define the board as a square matrix $A = (a_{i,j})$, $a_{i,j} = \pm 1$, with $1 \le i,j \le n$. $a_{i,j} = 1$ represents the *(i−j)−th* light on, i.e. the light on row *i*, column *j* is on, while $a_{i,j} = -1$ represents the *(i−j)−th* light off. The **column switches** are denoted by

*Variants of the Gale-Berlekamp Switching Game and their Solutions: Balancing the Rectangle and the Cube*

variables $x_1, x_2, x_3, \cdots, x_n = \pm 1$. $x_j = +1$ denotes that the $j-th$ column switch is not switched while $x_j = -1$ denotes that the $j-th$ column switch is switched once. Similarly, the **row switches** are denoted by variables $y_1, y_2, y_3, \cdots, y_n = \pm 1$.

## 2 AIMS AND OBJECTIVES

In this study, we aim to find out whether the algorithm of J. Beck and J. Spencer can be adapted to rectangle boards and cubes; then we would explore

- the specific conditions under which the algorithm can be applied to find the minimum imbalance in the specific boards; and
- the minimum imbalance of any given specific board satisfying these conditions.

## 3 METHODOLOGY

The proof of J. Beck and J. Spencer [1] provided a step-by-step algorithm which allows one to find not only the minimum imbalance in the board, but also the row and column switches to achieve this. Therefore, we hope to integrate the algorithm into Python and construct a complete Python routine which will be able to find the column and row switches to minimise the imbalance in the board and calculate the imbalance, for any given board.

We then used Python to randomly generate initial configurations of the board and test run our algorithm. The configuration of the randomly generated board and the imbalance were collected in pairs. This allowed us to analyse the conditions under which the algorithm works or fails, and subsequently sort out the special cases and solve them separately.

## 4 RESULTS AND DISCUSSION

### 4.1 Balancing the Rectangle

In the paper, the authors proved that the imbalance could be reduced to at most by turning switches for some rows and/or columns, for any $n \times n$ square matrix. Using similar methods, we managed to prove that for any rectangle board $m \times n$ matrix with $n$ even and $m \leq n$, there exists $x \in (\pm 1)^n$ and $y \in (\pm 1)^m$ such that $|yAx| \leq 2$. i.e. there exist column and row switches to reduce the imbalance in the board to at most $2$.

Let the row sums of a board/matrix A be $r_1, r_2, \cdots, r_m$, with $1 \leq i \leq m$. Note that $r_i$ could be positive, zero or negative. In our proof, we will relax the condition that the row sums are nonnegative. This could be achieved simply by flipping the rows whose sums are negative.

**Definition (Property I).** *A matrix A has* Property I *if* $0 \leq r_i \leq 2(i-1)$, $\forall\ 1 \leq i \leq m$.

**Lemma 1.** *Given any* $m \times n$ *board/matrix* A*, we can switch some rows and columns so that* A *has* Property I. *We can find the switches in* $O(n^4)$ *time*.

When $m < n$, we augment the matrix $A$ with more rows so that we obtain a square matrix. Clearly, if the augmented matrix has *Property I*, then the original matrix $A$ has *Property I*. Hence, we will prove for the case $m = n$. i.e.

**Lemma 1A.** *Given an* $n \times n$ *matrix* A*, then there exists* $x \in (\pm 1)^n$ *such that* $Ax = (r_1, r_2, \ldots, r_n)^T$ *with* $|r_i| \leq 2(i-1)$, $\forall\ 1 \leq i \leq n$.

Define the bounded cube $H_n = \{x: -1 \leq x_i \leq 1, \forall\ i\}$. For any $x \in H_n$, we say that a component $x_i$ is *floating* if $x_i \neq \pm 1$.

For $J \subseteq \{1, 2, \ldots, n\}$ and $1 \leq i \leq n$, we define $A_{i,J}$ to be the submatrix of $A$ whose rows are $1, 2, ..., i$ and columns belong to $J$. We use $A_i$ to denote the $i \times n$ matrix $A_{i,\{1, \ldots, n\}}$.

First, we set $x^{(0)} = 0$, and so all components of $x$ are floating. Set $J = \{1, 2, \ldots, n\}$. We consider $M = A_{n-1,J} = A_{n-1}$ and look at the subspace $U = \{c: Mc = 0\}$.

Since $M$ has more columns than rows, $U$ has dimension at least one. In other words, there exists a nonzero $c \in U$ and so, $A_{n-1}(\lambda c) = 0$ for all $\lambda \geq 0$. We choose the largest $\lambda$ so that $\lambda c$ belongs to $H_n$ and set $x^{(1)}$. Note that by the choice of $\lambda$, we have that $x^{(1)}$ has at least one less floating component.

Then we repeat this process to reduce the number of floating components. In the $i-th$ iteration, we will set $J$ to be the set of indices of the floating components of $x^{(i-1)}$. Set $M = A_{|J|-1,J}$ and consider the subspace $U = \{c: Mc = 0\}$.

As before, $U$ has dimension at least one and so, we can find a nonzero vector in $c$. Abusing notation, we insert $0$s into $c$ at positions not in $J$ so that $c$ has length $n$. Now, we check that $A_{|J|-1} (x^{(i-1)} + \lambda c) = 0$. Choose the largest $\lambda$ so that $x^{(i-1)} + \lambda c$ belongs to $H_n$ and set $x^{(i)} = x^{(i-1)} + \lambda c$. Note that by the choice of $\lambda$, we have that the number of floating components in $x^{(i)}$ is smaller than the number of floating components in $x^{(i-1)}$.

In the worst case, we perform $n-1$ iterations and $x^{(n-1)}$ has at most one floating component and $A_1 x^{(n-1)} = 0$. Now, since $n$ is even, it is not difficult to argue that the number of floating components must be zero.

We will prove that, with $x^{(n-1)}$, $Ax^{(n-1)}$ satisfies **Lemma 1A**. Indeed, for any $i$, during the procedure, as long as $|J| > i$, we have $A_{|J|-1} x = 0$, hence $r_i = 0$. Therefore, $r_i$ will only be changed if $|J| \leq i$. Through every iteration, all non-floating elements of $x$ will be preserved, and thus there is maximum $i$ elements of $x$ changed from the first time that $|J| \leq i$ until $|J| = 0$ to give the final $x^{(n-1)}$.

Now, for any $k$, if $x_k$ is floating, $|a_{i,k} x_k| < 1$, hence for each element of $x$ changed, $r_i$ changes by less than $2$. When a maximum of $i$ elements of $x$ is changed, $r_i$ changes by less than $2i$ i.e. $|r_i| < 2i$.

Since $n$ is even, $r_i$ is even, thus $|r_i| \leq 2(i-1)$, $\forall\ 1 \leq i \leq n$, satisfying **Lemma 1A**.

After finding the column shifts and switching the negative row sums so that all $r_i \geq 0$, we now proceed on to **Lemma 2**.

**Definition (Property II)**. *A sequence of* m *nonnegative numbers* $s_1, s_2, \ldots, s_m$ *has* Property II *if* $s_1 \leq 2$ *and* $s_i \leq s_1 + \cdots + s_{i-1} + 2$ *for all* $2 \leq i \leq m$.

**Lemma 2.** *Given any* m × n *board/matrix* A. *If* A *has* Property I, *then we can reorder the rows and switch one column and some rows so that the* m *row sums have* Property II. *We can find the switches in* $O(n^2)$ *time.*

First, we reorder the rows so that the row sums are nondecreasing. Note that after reordering, *Property I* still holds. In other words, $r_i \leq 2(i-1)$ for all $i$. Set $u$ and $v$ to be the number of row sums that are $0$ and $2$ respectively. In other words,

- $r_1 = \cdots = r_u = 0$
- $r_{u+1} = \cdots = r_{u+v} = 2$
- $r_i \geq 4$ for $i \geq u + v + 1$

There are two possible cases.

- If $v \geq n/2$, then by *Property I*, we can easily check that *Property II* is satisfied.

  - For $i \leq u + v$, we have $r_i \leq 2$ and so, less than or equal to the respective partial sum plus $2$.
  - For $i \geq u + v + 1$, we have that the partial sum is at least $2v \geq n$. Since any row sum is at most $n$, the condition is true.

- If $v < n/2$, there must be one column (say column $j$) such that in these $v$ rows, the number of $-1$s is less than or equal to the number of $+1$s in column $j$ (or else the total number of $-1$s in these $v$ rows would be greater than the total number of $+1$s, leading to at least one negative row sum i.e. contradiction). Let $r'_i$ be the absolute value of the new $r_i$ after column $j$ is flipped. We would have:

  - Since $r_1 = r_2 = r_3 = \cdots = r_u = 0$, after switching, $r'_1 = r'_2 = r'_3 = \cdots = r'_u = 2$.
  - $r'_{u+1} + r'_{u+2} + r'_{u+3} + \cdots + r'_{u+v} \geq r_{u+1} + r_{u+2} + r_{u+3} + \cdots + r_{u+v} = 2v$
  - For $i \geq u + v + 1$, since $4 \leq r_i \leq 2(i-1)$, we have $2 \leq r'_i \leq 2i$.

We prove that these $r'_i$ values satisfy *Property II*. Indeed:

- $r'_1 = 2$
- $r'_1 + r'_2 + r'_3 + \cdots + r'_i + 2 \geq 4 \geq r'_{i+1}$, for all $1 \leq i < u + v$
- $r'_1 + r'_2 + r'_3 + \cdots + r'_i + 2 \geq 2i + 2 \geq r'_{i+1}$ for all $i \geq u + v$

Thus **Lemma 2** is proven. We can now find the row switches.

**Lemma 3. (Greedy Technique)** *If the* $m$ *nonnegative numbers* $s_1, s_2, \ldots, s_m$ *have* Property II, *there is a linear-time method to find* $y_1, y_2, \ldots, y_m \in \{\pm 1\}$ *so that* $|y_1 s_1 + y_2 s_2 + \cdots + y_m s_m| \leq 2$.

For $1 \leq i \leq m$, set $S_i = s_1 + s_2 + \cdots + s_i$ and set $S_0 = 0$. Hence, *Property II* implies that $s_i \leq S_{i-1} + 2$ for all $1 \leq i \leq m$.

The **Greedy Technique** described by Beck and Spencer [1] does the following:

- Set $y_m = 1$ and $\sigma_m = y_m s_m = s_m$.

- Suppose that $y_i, y_{i+1}, \ldots, y_m$ are found and we have the values $\sigma_i, \sigma_{i+1}, \ldots, \sigma_m$.

    - If $\sigma_i \geq 0$, set $y_{i-1} = -1$ and $\sigma_{i-1} = \sigma_i - s_{i-1}$.
    - If $\sigma_i < 0$, set $y_{i-1} = 1$ and $\sigma_{i-1} = \sigma_i + s_{i-1}$.

Now, using backwards induction, we can show that statement **P**: $|\sigma_i| \leq S_{i-1} + 2$ is true for $1 \leq i \leq m$. Indeed:

- When $i = m$, $|\sigma_m| = s_m \leq S_{m-1} + 2$. **P** is true for $i = m$.

- Suppose **P** is true for $i = k$ for $2 \leq k \leq m$ i.e. $|\sigma_k| \leq S_{k-1} + 2$. We will prove that **P** is also true for $i = k-1$ i.e. $|\sigma_{k-1}| \leq S_{k-2} + 2$.

    - If $\sigma_k \geq 0$, $\sigma_{k-1} = \sigma_k - s_{k-1}$. Since $|\sigma_k| \leq S_{k-1} + 2$,

    $$0 \leq \sigma_k \leq S_{k-1} + 2$$
    $$\Rightarrow -s_{k-1} \leq \sigma_k - s_{k-1} \leq S_{k-1} - s_{k-1} + 2$$
    $$\Rightarrow -S_{k-2} - 2 \leq -s_{k-1} \leq \sigma_{k-1} \leq S_{k-2} + 2$$

    Therefore, $|\sigma_{k-1}| \leq S_{k-2} + 2$

    - Similarly, if $\sigma_k < 0$, we can also prove that $|\sigma_{k-1}| \leq S_{k-2} + 2$.

Hence, by mathematical induction, **P** is true for all $1 \leq i \leq m$. Therefore, at the end, $|\sigma_1| = |y_1 s_1 + \cdots + y_m s_m| \leq S_0 + 2 = 2$.

Applying the **Greedy Technique** and **Lemma 3** on $(r'_1, r'_2, \ldots, r'_m)$ which satisfies *Property II*, we can find $y \in (\pm 1)^m$ such that **imbalance** $= |y_1 r'_1 + \cdots + y_m r'_m| \leq 2$. This completes the proof.

### 4.2 Balancing the Cube

In the proof above, we balanced the rectangle by treating it as a two-dimensional array. At this point, it is natural to expand this problem to a three-dimensional array i.e. turning the board into a cube of lights. In this new configuration, instead of having row and column switches, we introduce $X$-, $Y$-, and $Z$- plane switches. Turning a switch would flip all the lights on the whole planar layer.

In this paper, we only consider the simplest three-dimensional configuration, a perfect $n \times n \times n$ cube. The objective of the game remains the same: to minimise the imbalance in the board. We define $P_n$ to be the minimum value that satisfies the following condition: for any initial state of an $n \times n \times n$ cube, there always exist $X$-, $Y$-, and $Z$- plane switches to reduce the imbalance in the cube to at most $P_n$.

Starting with the most basic case, a $n \times n \times n$ cube, we have the following theorem.

**Theorem 2**. For a $2 \times 2 \times 2$ cube, $P_2 = 2$.

First, we prove that for any initial state, there always exist $X$-, $Y$- and $Z$- plane switches to reduce the imbalance in the cube to at most $2$.

Fortunately, after some simplifications below, a $2 \times 2 \times 2$ cube is simple enough that we can manually verify every possible initial configuration to find the minimum imbalance.

Since there are $8$ lights, the minimum imbalance must always be an even integer. Hence we only need to consider the cases when the initial imbalance is $4$, $6$, or $8$. We can then easily calculate the number of lights on and off for each case, and this leaves us only 5 possible initial states to deal with (Figure 1).

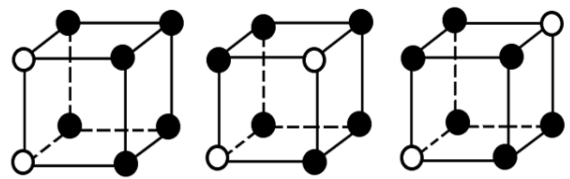

imbalance = 4

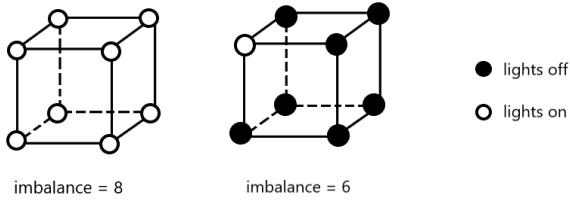

Figure 1. Unsolved initial states of a *2 x 2 x 2* cube

Each of them can easily be "solved" with a single switch (example of one case shown in Figure 2).

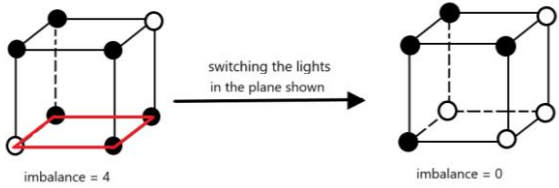

Figure 2. Example "solve" of one case

This means $P_2 = 2$. We now provide an example of a cube whose minimum imbalance is exactly *2* (Figure 3).

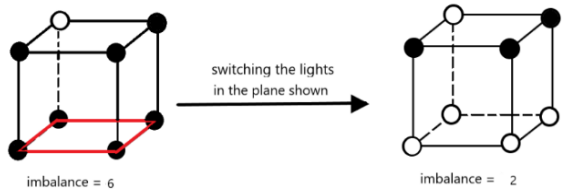

Figure 3. Example cube with minimum imbalance = *2*

Any plane switch would flip *4* lights, and the imbalance changes by *0*, *4*, or *8*. As the initial imbalance is *2*, the imbalance at any point in time is always congruent to *2* (modulo *4*), and cannot be reduced to *0*.

With this, we now move on to the more complex *4 x 4 x 4* cube.

**Theorem 3**. For a *4 x 4 x 4* cube, $P_4 = 4$.

We introduce the notations as follows.

We index the cube with the coordinates *(i, j, k)* where $0 \leq i, j, k \leq 3$. The state of the light at *(i, j, k)* is denoted as $A_{ijk}$. $A_{ijk} = 1$ means the light at *(i, j, k)* is on while $A_{ijk} = -1$ means the light at *(i, j, k)* is off.

For convenience, we will write this *4 x 4 x 4* cube as a *16 x 4* array (Figure 4).

Here, we have twelve switches $x_0, x_1, x_2, x_3, y_0, y_1, y_2, y_3, z_0, z_1, z_2, z_3$, that have the following effects on the cube (Figure 4).

- Switching $x_i$ changes $A_{ijk}$ for all *j,k* i.e. flips all the lights on the $X_i$ planar layer.

- Switching $y_j$ changes $A_{ijk}$ for all *k,i* i.e. flips all the lights on the $Y_j$ planar layer.

- Switching $z_k$ changes $A_{ijk}$ for all *i,j* i.e. flips all the lights on the $Z_k$ planar layer.

Hence, using $x_i, y_j, z_k = \pm 1$ to represent whether the switch $x_i, y_j, z_k$ respectively is thrown, then the imbalance is then given by $\sum_{i,j,k} x_i y_j z_k A_{ijk}$.

$$\begin{array}{cccc|cccc|cccc|cccc} A_{000} & A_{100} & A_{200} & A_{300} & A_{010} & A_{110} & A_{210} & A_{310} & A_{020} & A_{120} & A_{220} & A_{320} & A_{030} & A_{130} & A_{230} & A_{330} \\ A_{001} & A_{101} & A_{201} & A_{301} & A_{011} & A_{111} & A_{211} & A_{311} & A_{021} & A_{121} & A_{221} & A_{321} & A_{031} & A_{131} & A_{231} & A_{331} \\ A_{002} & A_{102} & A_{202} & A_{302} & A_{012} & A_{112} & A_{212} & A_{312} & A_{022} & A_{122} & A_{222} & A_{322} & A_{032} & A_{132} & A_{232} & A_{332} \\ A_{003} & A_{103} & A_{203} & A_{303} & A_{013} & A_{113} & A_{213} & A_{313} & A_{023} & A_{123} & A_{223} & A_{323} & A_{033} & A_{133} & A_{233} & A_{333} \end{array}$$

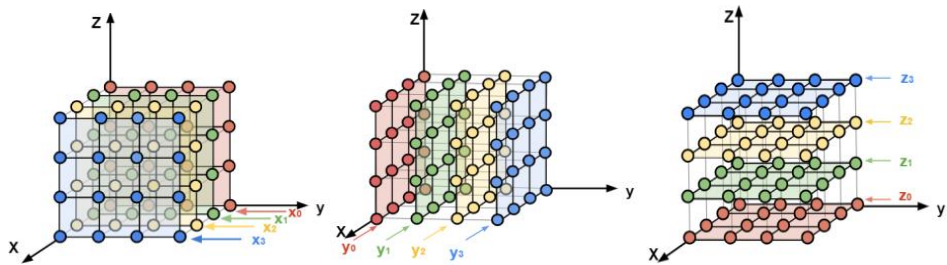

Figure 4. Configuration of a *4 x 4 x 4* cube

Define *X*- sums, *Y*- sums and *Z*- sums as follows.

$$X_i = \sum_{j,k} A_{ijk}$$

$$Y_j = \sum_{i,k} A_{ijk}$$

$$Z_k = \sum_{i,j} A_{ijk}$$

Note that $Z_k$ denotes the sum of the rows in our *16 x 4* array representation.

First, we will prove the following lemma.

**Lemma 4**. *For all initial states of a* 4 x 4 x 4 *cube, there always exist* X-, Y-, *and* Z- *plane switches to reduce the imbalance in the cube to at most* 4.

Without loss of generality, we're concerned with switching the *Z*- planes. We will only consider positive *Z*- sums, since we can always switch the corresponding *z*- switch if a *Z*- sum is negative. Also note that the exact order of the layers does not matter, and we can always reorder them if necessary.

**Lemma 4A**. *Given any initial state, we can throw some* x- *and* y- *switches so that*

$$Z_0 + Z_1 + Z_2 + Z_3 \leq 14.$$

We consider a *4 x 4* two-dimensional light board. For each possible combination of $x_i$ and $j_k$ switches, we calculate the imbalance in the board, and then we averaged this imbalance value over all $2^8$ possible combinations of $x_i$ and $j_k$ switches. We then repeat this process for all possible initial states of the light board, and find that the maximum average imbalance value for all possible initial states is *3.5*.

A cube contains four such layers, thus using an averaging argument, the maximum average imbalance value of any whole cube over all $2^8$ possible $x_i$ and $j_k$ switches is *3.5 x 4 = 14*.

Thus we can always find a combination of $x_i$ and $j_k$ switches such that the imbalance of the cube is no more than *14* i.e. $Z_0 + Z_1 + Z_2 + Z_3 \leq 14$. This completes the proof of **Lemma 4A**.

Recall and modifying **Lemma 3 (Greedy Technique)**, we have **Lemma 4B**.

**Lemma 4B. (Greedy Technique).** *If*

$Z_0 \leq 4$,

$Z_1 \leq Z_0 + 4$,

$Z_2 \leq Z_0 + Z_1 + 4$,

$Z_3 \leq Z_0 + Z_1 + Z_2 + 4$,

*then we can throw the* z- *switches so that the imbalance is at most* 4.

With **Lemma 4A** and checking all possible combinations, there are only *10* cases which we need to manually check for which the conditions of **Lemma 4B** does not hold and the **Greedy Technique** does not apply (Table 1).

Table 1. Cases to be manually checked

| Number | Case | No. of *-1*s in planar layer | | | |
|---|---|---|---|---|---|
| | | $Z_0$ | $Z_1$ | $Z_2$ | $Z_3$ |
| 1 | [0, 0, 0, 8] | 8 | 8 | 8 | 4 |
| 2 | [0, 0, 2, 8] | 8 | 8 | 7 | 4 |
| 3 | [0, 0, 0, 14] | 8 | 8 | 8 | 1 |
| 4 | [0, 0, 0, 12] | 8 | 8 | 8 | 2 |
| 5 | [0, 0, 2, 12] | 8 | 8 | 7 | 2 |
| 6 | [0, 0, 0, 6] | 8 | 8 | 8 | 5 |
| 7 | [0, 0, 0, 10] | 8 | 8 | 8 | 3 |
| 8 | [0, 0, 4, 10] | 8 | 8 | 6 | 3 |
| 9 | [0, 0, 2, 10] | 8 | 8 | 7 | 3 |
| 10 | [0, 2, 2, 10] | 8 | 7 | 7 | 3 |

One key observation we can make is the effect when an *x*- or *y*- switch is thrown (Table 2).

Table 2. Effect on $Z_i$ upon throwing one switch

| No. of *-1*s in row/column | Effect on $Z_i$ |
|---|---|
| 0 | -8 |
| 1 | -4 |
| 2 | 0 |
| 3 | +4 |
| 4 | +8 |

Further using the **Greedy Technique** and checking all possible combinations, we can then come up with **Lemma 4C** to simplify our analysis.

**Lemma 4C.** *Suppose that* $[Z_0, Z_1, Z_2, Z_3]$ *obeys the following conditions.*

- $Z_i \leq 8$ *for* $i = 0, 1, 2, 3$
- $Z_i \in \{2, 6\}$ *for at most two* $i$s
- $[Z_0, Z_1, Z_2, Z_3] \neq [0, 0, 0, 6], [0, 0, 0, 8], [0, 0, 2, 8], [0, 6, 8, 8], [0, 8, 8, 8], [2, 8, 8, 8]$

*Then we can throw z- switches to achieve an imbalance of at most 4.*

With **Lemma 4C**, our direction for solving these cases is to switch one or two *x*- or *y*- switch(es) so that $[Z_0, Z_1, Z_2, Z_3]$ satisfies **Lemma 4C**'s conditions.

For **Case 1** *[0, 0, 0, 8]* and **Case 2** *[0, 0, 2, 8]*, we solve by finding one single switch to reduce $Z_3$ to *4* (through finding one row or column with only one *-1* in the $Z_3$ planar layer), while ensuring $Z_0, Z_1, Z_2 \leq 8$. Consider two cases for the $Z_3$ planar layer.

- At least one row or column has only one *-1*.

By considering the number and configurations of *-1*s in each planar layer, we see that for each case, the minimum number of rows or columns with only one *-1* in the $Z_3$ plane (which is also the number of possible switches to reduce $Z_3$ to *4*) exceeds the maximum number of rows or columns in $Z_2$ that we must avoid switching to satisfy **Lemma 4C** (Table 3).

Table 3. **Case 1** and **Case 2** analysis

| Case | No. of $-1$s in $Z_3$ | Min no. of possible switches | Max no. of switches to avoid |
|---|---|---|---|
| 1 | 4 | 2 rows or 2 columns | 0 |
| 2 | 4 | 2 rows or 2 columns | 1 row and 1 column |

One example is provided for **Case 2** (Figure 5) where the possible switches and the switches to avoid are highlighted.

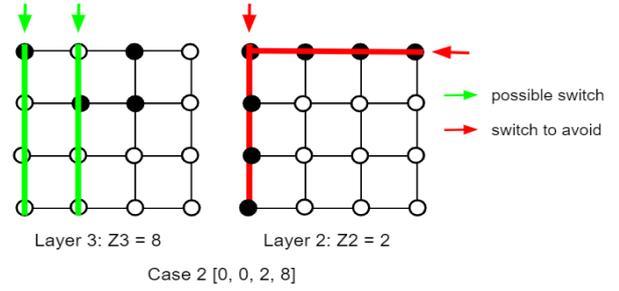

Figure 5. Example switches for **Case 2**

- There are no rows or columns with only one *-1*. The four *-1*s must be arranged in a square (Figure 6).

$$\begin{pmatrix} +1 & +1 & +1 & +1 \\ +1 & +1 & +1 & +1 \\ +1 & +1 & -1 & -1 \\ +1 & +1 & -1 & -1 \end{pmatrix}$$

Figure 6. Configuration of $Z_3$ planar layer

To reduce $Z_3$ to *4*, we must throw two switches as follows (Figure 7), while ensuring $Z_0, Z_1, Z_2 \leq 8$. Since it's not possible to turn $Z_i$ from 0 or 2 to 16 using one *x*- switch and one *y*- switch, we only need to avoid turning $Z_i$ to 10, 12 or 14. Examples of pairs of switches to avoid are highlighted in Figure 8.

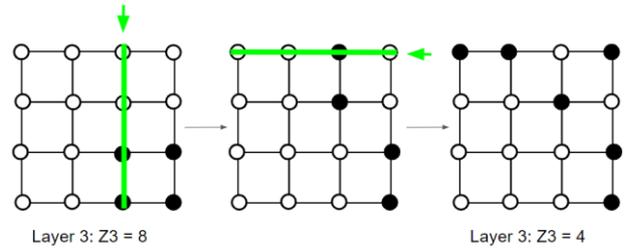

Figure 7. Example switching sequence to reduce $Z_3$ to *4*

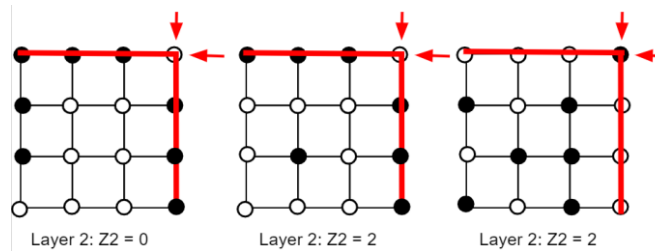

Figure 8. Pairs of switches to avoid for $Z_2$ planar layer

In every case, the minimum number of possible switches or switch pairs exceeds the maximum number of switches or switch pairs that must be avoided to satisfy **Lemma 4C** (Table 4). In other words, we can always switch one *x-* and *y-* switch so that the cube satisfies **Lemma 4C**, and can thus be balanced.

Table 4. **Case 1** and **Case 2** analysis (continued)

| Case | Min no. of possible pairs of switches | Max no. of pairs to avoid |
|---|---|---|
| 1 | 8 | 3 (1 for each planar layer) |
| 2 | 8 | 3 (1 for each planar layer) |

For **Case 3** *[0, 0, 0, 14]*, **Case 4** *[0, 0, 0, 12]*, **Case 5** *[0, 0, 2, 12]*, and **Case 6** *[0, 0, 0, 6]* we solve by finding one single switch to reduce $Z_3$ by *8* (or, in **Case 6**, by either *4* or *8*), while avoiding the undesirable cases in **Lemma 4C**. Specifically for **Case 6**, taking into account **Case 2** *[0, 0, 2, 8]* which has been solved, we just need to ensure $Z_2 \neq 8$.

Similar to the above cases, the minimum number of possible switches exceeds the maximum number of switches that must be avoided to satisfy **Lemma 4C** (Table 5). The cube can thus be balanced.

Table 5. **Case 3**, **Case 4, Case 5** and **Case 6** analysis

| Case | No. of −1s in $Z_3$ | Min no. of possible switches | Max no. of switches to avoid |
|---|---|---|---|
| 3 | 1 | 3 rows and 3 columns | 3 |
| 4 | 2 | 4 | 0 |
| 5 | 2 | 4 | 2 |
| 6 | 5 | 2 rows AND 2 columns | 2 rows OR 2 columns OR 1 row and 1 column |

One example is provided for **Case 6** (Figure 9) where the possible switches and the switches to avoid are highlighted.

Having dealt with **Case 1** *[0, 0, 0, 8]*, **Case 2** *[0, 0, 2, 8]*, and **Case 6** *[0, 0, 0, 6]*, we can improve on **Lemma 4C** as follows:

**Lemma 4D**. *Suppose that* [$Z_0$, $Z_1$, $Z_2$, $Z_3$] *obeys the following conditions*.

- $Z_i \leq 8$ *for* i = 0, 1, 2, 3
- $Z_i \in \{2, 6\}$ *for at most two* i*s*
- [$Z_0$, $Z_1$, $Z_2$, $Z_3$] ≠ [0, 6, 8, 8], [0, 8, 8, 8], [2, 8, 8, 8]

*Then we can throw z- switches to achieve an imbalance of at most* 4.

For **Case 7** *[0, 0, 0, 10]*, **Case 8** *[0, 0, 4, 10]*, and **Case 9** *[0, 0, 2, 10]*, we solve by finding one single switch to reduce $Z_3$ by either *4* or *8*, while avoiding the undesirable cases in **Lemma 4D** by ensuring $Z_0$, $Z_1$, $Z_2 \leq 8$ and at most one of $Z_0$, $Z_1$, $Z_2 = 8$.

The minimum number of possible switches exceeds the maximum number of switches that must be avoided to satisfy **Lemma 4D** (Table 6). The cube can thus be balanced.

Table 6. **Case 7, Case 8,** and **Case 9** analysis

| Case | No. of −1s in $Z_3$ | Min no. of possible switches | Max no. of switches to avoid |
|---|---|---|---|
| 7 | 3 | 3 rows and 3 columns | 3 (similar to **Case 3**) |
| 8 | 3 | 3 rows and 3 columns | 5 (2 for $Z_2 \leq 4$ and 3 for $Z_1 \leq 4$) |
| 9 | 3 | 3 rows and 3 columns | 2 (for $Z_2 \neq 10$) |

One example is provided for **Case 8** (Figure 10) where the possible switches and the switches to avoid are highlighted.

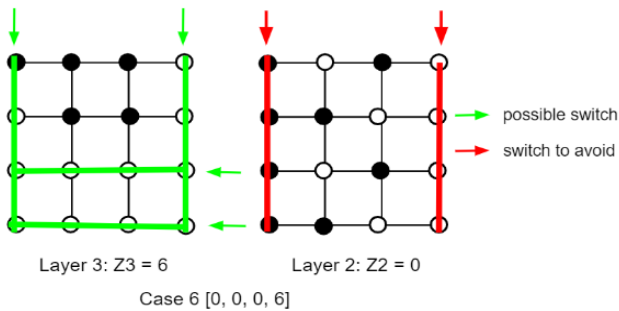

Figure 9. Example switches for **Case 6**

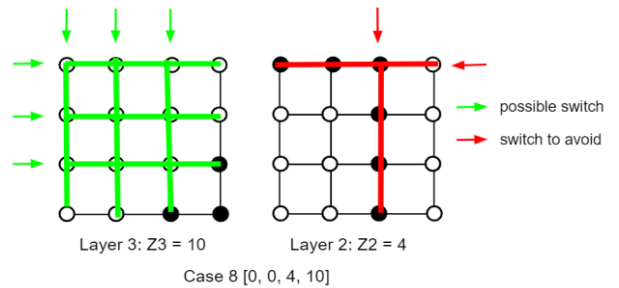

Figure 10. Example switches for **Case 8**

For **Case 10** *[0, 2, 2, 10]*, **Lemma 4D** does not apply, however we can manually check every case and expand **Lemma 4D** as follows: the cube can always be balanced if $Z_i$ satisfies

- $Z_0 \in \{0, 4, 8\}$
- $Z_1, Z_2, Z_3 \in \{2, 6\}$
- $[Z_0, Z_1, Z_2, Z_3] \neq [0, 6, 6, 6]$

We solve this case by finding one switch to reduce $Z_3$ by either *4* or *8*, while ensuring $Z_1, Z_2 < 10$ (by avoid switching a row or column full of *-1*s on the corresponding planar layer) and $Z_2 \neq 6$ to avoid the undesirable case *[0, 6, 6, 6]*.

Without loss of generality, assume that the total number of rows or columns that are full of *-1*s in the $Z_2$ planar layer is greater than or equal to that in $Z_1$.. We consider three cases for the number of such rows or columns on the $Z_2$ planar layer. In all three cases, the minimum number of possible switches exceeds the maximum number of switches that must be avoided (Table 7).

For each sub-case, the "worst case scenario" i.e. the initial configuration of the $Z_2$ planar layer which results in the maximum number of switches to avoid is illustrated in Figure 11, with the switches to be avoided highlighted in red.

Thus for **Case 10**, the cube can always be balanced.

Having checked every case, we complete the proof for **Lemma 4**, and from this we can deduce that $P_4 \leq 4$.

We now provide an example of a cube whose minimum imbalance is exactly *4* (Figure 12).

To prove the minimum imbalance, we introduce the following Property and Lemma.

**Definition (Property III).** *A cube has* Property III *if $x_i, y_j, z_k$ are all divisible by* 4 *for all* $0 \leq i, j, k \leq 3$.

**Lemma 5.** *After any number of switches, the imbalance of a cube with* Property III *will always give the same remainder when divided by* 8.

Table 7. **Case 10** analysis

| Case | No. of $-1$s in $Z_3$ | Min no. of possible switches | No. of rows/columns full of $-1$s in $Z_2$ | Max no. of switches to avoid |
|---|---|---|---|---|
| 10a | 3 | 3 rows and 3 columns | 1 row AND 1 column | 4 (2 for $Z_2 < 10$, 2 for $Z_1 < 10$, 0 for $Z_2 \neq 6$) |
| 10b | 3 | 3 rows and 3 columns | 1 row OR 1 column | 4 (1 for $Z_2 < 10$, 1 for $Z_1 < 10$, 2 for $Z_2 \neq 6$) |
| 10c | 3 | 3 rows and 3 columns | 0 | 5 (0 for $Z_1, Z_2 < 10$, 5 for $Z_2 \neq 6$) |

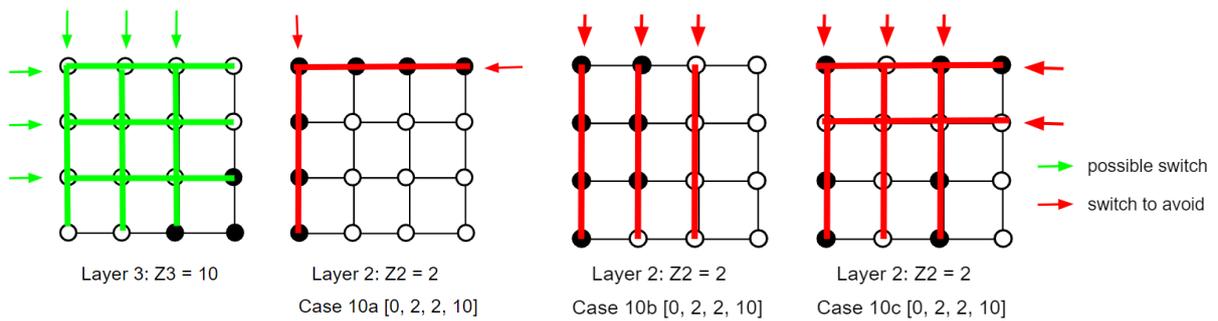

Figure 11. Maximum number of switches to avoid for each sub-case of **Case 10**

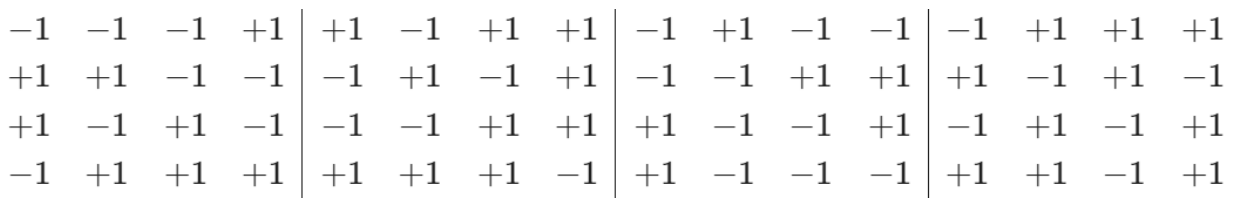

Figure 12. Example *4 x 4 x 4* cube whose minimum imbalance = *4*

We consider the effect of one switch on the cube. Without loss of generality, assume a $Z_k$ switch was switched, turning the plane sum $Z_k$ into $-Z_k$. The change in imbalance of the board is $-2Z_k$, which is divisible by *8*, due to *Property III*. The imbalance thus retains its remainder when divided by *8*.

We now prove that after one switch, *Property III* is retained. The absolute value of *Z-* sums remain unchanged, thus all new *Z-* sums are divisible by *4*.

Consider *X-* sums. When a $z_k$ switch is switched, in each $X_i$ layer, *4* lights (in a *4* x *1* x *1* block) are flipped. Each *X-* sum changes by *0*, *4*, or *8*. Since all old $X_i$s are divisible by *4*, all new $X_i$s are also divisible by *4*. The same argument can then be applied to *Y-* sums. Thus the new cube also satisfies *Property III*.

Repeating the same argument on the new cube, and any subsequent cubes, we eventually reach the conclusion that after any number of switches, the cube retains *Property III*, and the imbalance always retains its remainder when divided by *4*. This completes **Lemma 5**'s proof.

Back to the example cube, we can verify that it satisfies *Property III*. The initial imbalance is *4*. By **Lemma 5**, the imbalance is always congruent to *4* (modulo *8)*, and thus the minimum achievable imbalance is exactly *4*. This ensures that $P_4 = 4$, and completes **Theorem 3**.

## 5 CONCLUSION

In this paper, we have proved the existence of row and column switches to reduce the imbalance of any *m x n* rectangular board, with $m \leq n$ and *n* even. This is achieved through a algorithm consisting of three distinct subroutines, that can all be easily be implemented with at most $O(n^4)$ complexity. We then expanded the game to a *n x n x n* cube with plane switches. The *2 x 2 x 2* cube is simple enough that its cases can be manually and exhaustively checked, and thus it can be verified that $P_2 = 2$. For the more complex *4 x 4 x 4* cube, we used an averaging argument to limit the plane sums. We then combined it with a modified version of the **Greedy Technique** to narrow down to *10* outstanding cases, which we manually checked to verify that $P_4 = 4$.

However, we have yet to explore the case for a rectangular *m x n* board, where $m \leq n$ and *n* is odd, or subsequently, when both *m* and *n* are odd, due to the added complexity of having an odd number of columns. This necessitates a closer inspection and refinement of the three subroutines. Moreover, for the three-dimensional cube, we have only examined the simplest configurations i.e. *n = 2* and *n = 4*. The minimum imbalance and general solution for the general *n x n x n* cube are undiscovered and remain a subject for future studies. We have strong belief that the proofs and techniques outlined in this paper would provide some valuable insights and inspirations that can be applied to the general *n x n x n* cube.

## ACKNOWLEDGEMENTS

This paper and the research behind it would not have been possible without the exceptional support of our Supervisor, Asst. Professor Kiah Han Mao. He provided us with the backbone for our Python routine, and his contributions towards the earlier draft of the report have greatly benefitted us in articulating the proofs in a clear and concise manner. We would like to show our deepest gratitude to Professor Kiah for his expertise, guidance and mentorship throughout the course of this research journey.

We would also like to thank our teacher mentors, Mr Teo Yik Tee and Mrs Lim Siew Eng, for their support throughout the project, and for ensuring that all the deadlines and criteria of the project are met. Without their help and instruction, this report would never have been possible.